\documentclass[11pt]{article}
\usepackage[latin1]{inputenc}
\usepackage{amsmath}
\usepackage{amsfonts}
\usepackage{amsthm}
\numberwithin{equation}{section}
\newtheorem{teo}{Theorem}[section]
\newtheorem{lem}[teo]{Lemma}

\newtheorem{cor}[teo]{Corollary}
\theoremstyle{definition}
\newtheorem{defi}{Definition}
\theoremstyle{remark}
\newtheorem{rem}{Remark}

\newcommand{\vv}{{\mbox{\boldmath$v$}}}
\newcommand{\ve}{\varepsilon}

\newcommand{\pical}{\mathcal{P}}
\newcommand{\lcal}{\mathcal{L}}

\newcommand{\R}{\mathbb{R}}
\newcommand{\Q}{\mathbb{Q}}
\newcommand{\F}{\mathfrak{F}}
\newcommand{\tto}{\rightarrow}

\newcommand\io{\int_{\Omega}}
\DeclareMathOperator{\sign}{sign}
\DeclareMathOperator{\Lip}{Lip}

\author{Luigi Ambrosio and Filippo Santambrogio\thanks{Scuola Normale Superiore, Piazza dei Cavalieri 7, 56126 Pisa, Italy, {\tt l.ambrosio@sns.it}, {\tt f.santambrogio@sns.it}}}
\title{Necessary optimality conditions for geodesics \\ in weighted Wasserstein spaces}
\date{}
\begin{document}
\maketitle
{\bf Abstract:} The geodesic problem in Wasserstein spaces with a metric perturbed by a conformal factor is considered, and necessary optimality conditions are estabilished in a case where this conformal factor favours the spreading of the probability measure along the curve. These conditions have the form of a system of PDEs of the kind of the compressible Euler equations. Moreover, self-similar solutions to this system are discussed.
\sloppy

\section{Introduction}

Let us consider a closed convex set $\Omega\subset\R^d$ and the set
$\pical(\Omega)$ of probability measures in $\Omega$. Given $p\in (1,\infty)$,
we denote by $W_p(\Omega)$ the subspace of measures with finite $p$-th moments, i.e.
$$
W_p(\Omega):=\left\{\mu\in\pical(\Omega)\,:\, \int_\Omega|x|^p\,d\mu<\infty\right\}.
$$
We endow $W_p(\Omega)$ with the canonical Wasserstein distance $W_p(\mu,\nu)$
of order $p$ (see \cite{AmGiSa}, \cite{Villani} for the basic facts about $W_p$).

It is well known that $W_p(\Omega)$ is a length space, and that (constant speed) 
geodesics of $W_p(\Omega)$ are in one to one correspondence with optimal transport plans,
via McCann's linear interpolation procedure (see for instance Proposition~7.2.2 of
\cite{AmGiSa}). 
Here we consider, instead, the case when the Wasserstein metric is perturbed
by a conformal factor $L(\mu)$: by minimizing
\begin{equation}\label{min1}
\int_0^1 L^2(\mu_t)|\mu'|^2(t)\,dt
\end{equation}
among all curves $\mu$ connecting $\mu_0=\mu$ to $\mu_1=\nu$, one obtains
a new squared distance depending on $p$ and $L$, and we are interested in computing the 
geodesics relative to this distance. In \eqref{min1}, $|\mu'|(t)$ is the rate of change of $W_p$
also called metric derivative, along the curve $\mu$, see \eqref{min2}.

This problem has been introduced in \cite{BraButSan}, where the main goal was to choose a factor $L$ 
favouring atomic measures in order to give a time-dependent approach to some branched transport problems 
which may be applied to the study of river networks, pipe systems, blood vessels, tree structures\dots . 
In fact, by setting $L(\mu)=\sum_i a_i^r$ (for $0<r<1$) if $\mu=\sum_i a_i\delta_{x_i}$ and $L=+\infty$ 
on measures which are not purely atomic, there is a strong link between this variational problem and those 
which were first presented in \cite{xia1} and \cite{MadMorSol} (the latter uses in fact a time-dependent 
approach, but by means of measures on the space of paths instead of paths on the space of measures). 
This choice of $L$ is in fact a local functional on measures which, among probability measures, favours 
the most concentrated ones. In \cite{BraButSan}, as a natural counterpart, the case of local functionals 
$L$ which prefer spread measures is considered as well and the two problems sound somehow specular. 
The aim of the present paper is in fact to consider this second problem and to find out optimality 
conditions in the form of PDEs.

In particular, we study in detail the case when $L(\mu)$ is the $\gamma$-th power of the
$L^q$ norm of the density of $\mu$ with respect to Lebesgue measure ${\cal L}^d$,
with $q>1$ and $\gamma>0$ given, and $L(\mu)=+\infty$ if $\mu$ is a singular
measure. So, geodesics with respect 
to the new metric tend to spread the density as much as possible. Denoting by
$u_t$ the density of $\mu_t$, we find that a necessary optimality condition for
geodesics is
(for $p=2$, see \eqref{system} for general $p$) 
\begin{equation}\label{compressible}
\frac{d}{dt}\left(K(t)\vv u\right)
+K(t)\nabla\cdot\left(\vv\otimes\vv u\right)+H(t)\nabla u^q=0,
\end{equation}
where $\vv_t$ is the tangent velocity field of $\mu_t$, linked to $u_t$ via the
continuity equation $\frac{d}{dt}u_t+\nabla\cdot (\vv_t u_t)=0$. Here $H(t)<0$ and
$K(t)>0$ are suitable functions depending only on the metric derivative
of $\mu_t$ and on $L(\mu_t)$. As Brenier pointed to us, this equation is
very similar to the compressible Euler equation, but with a negative
pressure field $p=H(t)u^q$; a similar equation, with $H$ constant and
$q=3$, recently appeared also in \cite{guionnet}, in the one-dimensional case. 
In fact the main difference appears in the relationship between the $L$ part and the speed part: 
here it is multiplicative, while in \cite{guionnet} it is additive, as we will explain in a while.

The appearence of the Euler equation as an optimality condition is not
very surprising, taking into account the approach developed,
in the incompressible case, by Brenier (first in a purely Lagrangian framework 
in \cite{Least}, \cite{Br1}, and then in a mixed Eulerian-Lagrangian one in 
\cite{Br2}, \cite{Br3}). In this connection, we mention that our derivation
of the optimality condition differs from \cite{Br1}, \cite{Br3}, where duality
is used to perform first variations, and uses instead a perturbation argument
directly at the level of the primal problem.

Due to the non-convex nature of this problem, we don't know of any sufficient
minimality condition for the geodesics. In this connection, one may notice that,
in the case $\gamma=q/2$ and $p=2$, we have
$$\inf_{\delta>0}\,\delta\!\int_\Omega u^q\,dx+\frac{1}{\delta}\int_\Omega|\vv|^2u\,dx
=2L(u{\cal L}^d)\left(\int_\Omega|\vv|^2 u\,dx\right)^{1/2}$$
and the minimal $L^2(\mu)$ norm of $\vv$ is strictly linked to the metric derivative.
This suggests a connection between the ``multiplicative'' model studied here and
in \cite{BraButSan}, and the ``additive'' model
$$
\min\left\{\int_0^1\int_\Omega u^q+|\vv|^2u\,dxdt:\ \frac{d}{dt}u+\nabla\cdot (\vv u)=0\right\}
$$
subject to Dirichlet conditions at $t=0$ and $t=1$. This additive model, in the case $q=3$, 
is exactly the one studied in \cite{guionnet} (in this connection, see also \cite{loeper}). 
Notice that this problem is {\it convex} in the pair $(u,\vv u)$.
It turns out, indeed, that the (necessary and sufficient, by the convex nature 
of the problem) optimality conditions for the additive model are very similar 
to \eqref{compressible}, the only difference being that $H$ and $K$ do not depend on time. 

In the last part of the paper we compute and characterize particular self-similar or
solutions of \eqref{compressible}.

\medskip
\noindent
{\bf Acknowledgements.} We warmly thank Y. Brenier for many useful comments on an earlier
version of this paper and G. Buttazzo for the interest towards this work, which in particular 
lead us to study the self similar solutions of Section 3.

\section{Optimality Conditions for Weighted\\ Wasserstein Geodesics}

\subsection{A new velocity vector field}

\begin{defi}
If we are given a Lipschitz curve $\mu:[0,1]\tto W_p(\Omega)$, we define velocity field of the curve 
any vector field $\vv:[0,1]\times\Omega\tto \R^d$ such that for a.e. $t\in[0,1]$ the vector field 
$\vv_t=\vv(t,\cdot)$ belongs to $[L^p(\mu_t)]^d$ and the continuity equation 
$$\frac{d}{dt}\mu_t+\nabla\cdot (\vv\mu_t)=0$$
is satisfied in the sense of distributions: this means that for all 
$\phi\in C^1_c(\Omega)$ and any $t_1<t_2\in[0,1]$ it holds
$$\int \phi\, d\mu_{t_2}-\int\phi \,d\mu_{t_1}=\int_{t_1}^{t_2}ds\io \nabla\phi\cdot \vv_s\, d\mu_s,$$
or, equivalently, in differential form:
$$\frac{d}{dt}\int \phi \,d\mu_{t}=\io \nabla\phi\cdot \vv_t\, d\mu_t\qquad\mbox{ for a.e. }t\in[0,1].$$
We say that $\vv$ is the \emph{tangent} field to the curve $\mu_t$ if, for a.e. $t$, 
$\vv_t$ has minimal $[L^p(\mu_t)]^d$ norm for any $t$ among all the velocity fields.
\end{defi}

It is now well known (see for instance Theorem~8.3.1 and Proposition~8.4.5 in \cite{AmGiSa})
that for any Lipschitz or absolutely continuous curve $\mu_t$ with values in $W_p(\Omega)$ 
there exists a unique tangent field and moreover it is characterized by
\begin{equation}\label{min2}
\Vert \vv_t\Vert_{L^p(\mu_t)}
=|\mu'|(t)=\lim_{h\tto 0}\frac{W_p(\mu_{t+h},\mu_t)}{|h|}
\quad\mbox{ for a.e. }t\in[0,1].
\end{equation}
The right hand side, in the equality above is the rate of change of $W_p$ along
the curve $\mu_t$, also called metric derivative of $\mu_t$.

We want now to investigate how velocity fields change if we modify the curve $\mu_t$.

\begin{teo}\label{newfield}
Let a Lipschitz function $\mu_t:[0,1]\tto W_p(\Omega)$
and a smooth function $T:[0,1]\times\Omega\tto\Omega$ be given, such that for any $t$ 
the function $T_t:=T(t,\cdot)$ is a diffeomorphism. Let us consider the new curve $\mu_t'$ 
given by $\mu'_t=(T_t)_{\sharp}\mu_t$. If $\vv_t$ is a velocity field for $\mu_t$, then the vector 
field $\vv'$ defined by 
$$\vv'_t\cdot\mu'_t=
(T_t)_{\sharp}\left[\left(\nabla T_t\cdot \vv_t+\frac{\partial T}{\partial t}\right)\mu_t\right]$$
is a velocity field for $\mu_t'$.
\end{teo}

\begin{proof}
We have
\begin{multline*}
\io \phi\, d\mu'_{t+h}-\io\phi\, d\mu'_t=\io \phi\circ T_{t+h} d\mu_{t+h} - \io \phi\circ T_{t}\, d\mu_{t}\\
=\io \left(\phi\circ T_{t+h} - \phi\circ T_{t}\right)d\mu_{t+h} + \io \phi\circ T_{t}\, d(\mu_{t+h}-\mu_t)\\
=\io \left(\int_t^{t+h}(\nabla\phi)\circ T_s \cdot \frac{\partial T}{\partial t}|_s \,ds\right)d\mu_{t+h}+
\int_t^{t+h}ds\io(\nabla\phi)\circ T_t \cdot \nabla T_t\cdot \vv_s\,d\mu_s,
\end{multline*}
where in the last equality we have used the fact that $\vv_t$ is a velocity field for $\mu$, with test function 
$\phi\circ T_{t}$. It is now convenient to divide by $h$, rewrite and pass to the limit as $h\tto 0$:
\begin{multline}\label{dopodividereperh}
\frac{\io \phi d\mu'_{t+h}-\io\phi d\mu'_t}{h}=\io d\mu_{t+h}\frac{1}{h}\int_t^{t+h}(\nabla\phi)
\circ T_s \cdot \frac{\partial T}{\partial t}|_s \,ds \\+\io(\nabla\phi)\circ T_t \cdot \nabla T_t\cdot \vv_t\,d\mu_t+
\frac{1}{h}\int_t^{t+h}ds\io\nabla\psi_t\cdot \left(\vv_s d\mu_s-\vv_td\mu_t\right),
\end{multline}
where $\psi_t=\phi\circ T_t$. In the first term on the right hand side the measures $\mu_{t+h}$ weakly converge to $\mu_t$, 
since $t\mapsto\mu_t$ is Lipschitz continuous, while the 
integrand uniformly converges as a function of the space variable $x$ to 
$(\nabla\phi)\circ T_t \cdot \frac{\partial T}{\partial t}$ as $h\tto 0$. Hence we get convergence of the integral. 
If we prove that the last term tends to zero at least for a.e. $t\in[0,1]$ we get the thesis, since then we would have
\begin{multline*}
\lim_{h\tto 0}\frac{\io \phi d\mu'_{t+h}-\io\phi d\mu'_t}{h}\\
=\io\left((\nabla\phi)\circ T_t \cdot \frac{\partial T}{\partial t}+(\nabla\phi)\circ 
T_t \cdot \nabla T_t\cdot\vv_t\right)\,d\mu_t=\io\nabla\phi\cdot\vv'_td\mu'_t,
\end{multline*}
and this is nothing but the differential version of the continuity equation for $\vv'$ and $\mu'$ 
(it remains to prove $\vv_t'\in L^p(\mu'_t)$ but this is straightforward since $T_t$ is a diffeomorphism and 
this allows to write down the densities and estimate them).
To prove that the last term vanishes at the limit we see that, for fixed $\psi\in \Lip(\Omega)$ the function
$$s\mapsto g_{\psi}(s):=\io\nabla\psi\cdot \vv_sd\mu_s=\frac{d}{ds}\io\psi\,d\mu_s$$
is $L^{\infty}$ since $\mu_t$ is a Lipschitz curve in $W_p(\Omega)$ and hence almost any $s\in[0,1]$ is a Lebesgue point. 
This allows to fix a negligible set $N\subset[0,1]$ such that any point $t\in[0,1]\setminus N$ is a Lebesgue point 
for all the functions $g_{\psi_{t_1}}$ for $t_1\in\Q$. We fix now $t\in[0,1]\setminus N$ and try to prove that the 
last integral in \eqref{dopodividereperh} tends to zero. For $t_1\in\Q$ it holds
\begin{multline*}
\left|\frac{1}{h}\int_t^{t+h}ds\io\nabla\psi_t\cdot \left(\vv_s d\mu_s-\vv_td\mu_t\right)\right|\\
\leq
\frac{1}{h}\int_t^{t+h}ds\left|\io\nabla(\psi_t-\psi_{t_1})\cdot \vv_s d\mu_s\right|
+\left|\io\nabla(\psi_t-\psi_{t_1})\cdot \vv_t d\mu_t\right|+\\
\left|\frac{1}{h}\int_t^{t+h}ds\io\nabla\psi_{t_1}\cdot \left(\vv_s d\mu_s-\vv_td\mu_t\right)\right|\\
\leq \Lip(\psi_t-\psi_{t_1}) \Lip_{W_p}(\mu)+\left|\frac{1}{h}\int_t^{t+h}ds\io\nabla\psi_{t_1}\cdot 
\left(\vv_s d\mu_s-\vv_td\mu_t\right)\right|.
\end{multline*}
In the last sum the second term tends to zero by the fact that $t$ is a Lebesgue point for $g_{\psi_{t_1}}$ 
and the first term may be made as small as we want by choosing $t_1$ close to $t$, since $\psi_s=\phi\circ T_s$ 
and both $\phi$ and $T$ are regular.
\end{proof}

\subsection{Derivation of the optimality conditions}

We consider the minimization problem presented in \cite{BraButSan}, i.e. finding a curve of measures in $W_p(\Omega)$ of minimal length according to a metric which, roughly speaking is the Wasserstein (infinitesimal) metric multiplied
by a conformal factor. Precisely, if we define for $q>1$ the functional
$$L_q(\nu)=\begin{cases}\io u^q d\lcal^d&\mbox{ if }\nu=u\cdot\lcal^d\\
                      +\infty          & \mbox{otherwise,}\end{cases}$$
we want to minimize
$$\int_0^1 L_q(\mu_t)|\mu'|(t)\,dt,$$
where $|\mu'|(t)$ is the metric derivative of the curve $\mu$ 
and the minimization occurs among all the Wasserstein-Lipschitz curves $t\mapsto\mu_t$ with given initial and 
final points,  i.e. $\mu_0$ and $\mu_1$ are given probability measures in $W_p(\Omega)$. We will always consider only the non trivial case $\mu_0\neq\mu_1$.
If we define $V(\mu,t)=\io |\vv_t|^pd\mu_t$, where $\vv$ is the tangent field to the curve $\mu_t$, 
we know that $|\mu'|(t)=V(\mu,t)^{1/p}$. We may generalize the functional we want to minimize by considering
$$\F(\mu):=\int_0^1 L_q(\mu_t)^{\alpha}\,V(\mu,t)^{\beta}dt$$
which reduces to the case studied in \cite{BraButSan} if $\alpha=1$ and $\beta=1/p$. 
Notice that in this case the functional does not change under reparametrization of curves, while if $\beta>1/p$ 
the minimization selects a particular parametrization. For $\beta\leq 1/p$ the existence of a minimum is not ensured. 
Anyway we do not deal here with existence results (see \cite{BraButSan}), but we only look for necessary optimality conditions.
We will consider variations of $\mu$ of the form 
$$\mu^{\ve}_t=(T^{\ve}_t)_{\sharp}\mu_t\quad\mbox{ with }
\quad T^{\ve}(t,x)=x+\ve\xi(t,x),\,\,\,T^{\ve}=id+\ve\xi(t,\cdot),$$
for arbitrary regular functions $\xi\in C^{\infty}_c([0,1]\times\Omega;\R^d)$.
In the end optimality conditions will be expressed through a system of PDEs: we will obtain the result 
after collecting some lemmas. What we want to do now is exploiting the fact that for a minimizing curve $\mu$ 
the following quantity must be minimal for $\ve=0$:
$$\F(\mu^{\ve}_t)=\left(\int_0^1 F_{\ve}(t)^{\alpha}V_{\ve}(t)^{\beta}dt\right),$$
provided we define $F_{\ve}(t)=L_q(\mu^{\ve}_t)$ and $V_{\ve}(t)=V(\mu^{\ve},t)$.
Since it is not completely easy to deal with the term $V_{\ve}(t)$, we will replace it by 
$\tilde{V_{\ve}}(t)$, with $\tilde{V_{\ve}}(t)$ given by
$$\tilde{V_{\ve}}(t)=\io |(\vv^{\ve})_t|^p\,d\mu^{\ve}_t.$$
Here the vector field $\vv^{\ve}$ is the one we get by Theorem \ref{newfield} when the map $T$ is 
given by $T^{\ve}$ and the initial field $\vv_t$ is the tangent field to $\mu_t$. In this way we have 
$\tilde{V_{\ve}}(t)\geq V_{\ve}(t)$ (since $\vv^{\ve}_t$ is not necessarily of minimal $L^p$ norm) 
but $\tilde{V_{\ve}}(0) = V_{\ve}(0)$. Thus we may switch to considering $\tilde{V_{\ve}}(t)$ 
instead of $V_{\ve}(t)$, getting
$$\tilde{\F}(\mu_t^{\ve})=\left(\int_0^1 F_{\ve}(t)^{\alpha}\tilde{V_{\ve}}(t)^{\beta}dt\right).$$
We will compute the derivative of $\tilde{\F}(\mu_t^{\ve})$ with respect to $\ve$ and get the 
conditions we are looking for.
\begin{lem}\label{deriF}
If $\mu$ is a curve given by $\mu_t=u_t\lcal^d$ and such that $\F(\mu)<+\infty$, then for almost any $t\in[0,1]$ it holds
$$\frac{d}{d\ve}F_{\ve}(t)=(1-q)\io (JT^{\ve}_t)'\left(\frac{u_t}{JT^{\ve}_t}\right)^qd\lcal^d.$$
In particular, if we compute the derivative at $\ve=0$, we have
$$\frac{d}{d\ve}F_{\ve}(t)|_{\ve=0}=(1-q)\io (\nabla\cdot\xi)u_t^qd\lcal^d.$$
Moreover, for $\ve$ sufficiently small (depending on $T$, but not on $t$)
the following inequality holds: 
$$\frac{d}{d\ve}F_{\ve}(t)\leq C L_q(\mu_t).$$
\end{lem}
\begin{proof}
We look at the integrand function in the definition of $F_{\ve}$: to do this it is necessary to look 
at the density of the measure $\mu^{\ve}_t$. Thanks to the change of variables formula,
this density can be easily seen to be given by
$$u^{\ve}_t=\frac{u_t}{JT^{\ve}_t}\circ (T^{\ve}_t)^{-1},$$
where $J$ stands for the Jacobian (this formula is a consequence of $T^{\ve}_t$ 
being a diffeomorphism at least for small $\ve$).
Thus, after changing variables, we have
$$F_{\ve}(t)=L_q(\mu^\ve_t)=\io \left(\frac{u_t}{JT^{\ve}_t}\right)^qJT^{\ve}_t\,d\lcal^d.$$
The derivative of the integral is given by
$$(1-q)(JT^{\ve}_t)'\left(\frac{u_t}{JT^{\ve}_t}\right)^q,$$
where $(JT^{\ve}_t)'$ stands for the derivative w.r.t. $\ve$ of $JT^{\ve}_t$.
This quantity may be easily estimated by $Cu_t^q$, since $1-a\leq JT^{\ve}_t\leq 1+a$ and $(JT^{\ve}_t)'\leq B$
for suitable constants $a$ and $B$. 
Since for almost any $t$ the function $u_t$ must belong to $L^q$ (because the functional we are minimizing is finite) 
we can apply the dominated convergence theorem and get the thesis. To obtain the derivative at $\ve=0$ 
it is sufficient to notice that $(JT^{\ve}_t)'|_{\ve=0}=\nabla\cdot\xi,$
which is well-known. The same estimate we used to get dominated convergence may be used to get 
the last inequality.
\end{proof}
In the next lemma we consider the term $\tilde{V_{\ve}}$. 

\begin{lem}\label{deriV}
If $\mu$ is a curve such that $\F(\mu)<+\infty$, then for almost any $t\in[0,1]$ it holds
\begin{equation}\label{deriVeq}
\frac{d}{d\ve}\tilde{V_{\ve}}(t)=p\!
\io \left|\nabla T^{\ve}_t\cdot\vv_t+\frac{\partial T^{\ve}}{\partial t}\right|^{p-2}\!\left(\nabla T^{\ve}_t\cdot\vv_t+
\frac{\partial T^{\ve}}{\partial t}\right)\cdot\left(\nabla\xi\cdot\vv_t+\frac{\partial \xi}{\partial t}\right)d\mu_t.
\end{equation}
In particular, if we compute the derivative at $\ve=0$, we have
$$\frac{d}{d\ve}\tilde{V_{\ve}}(t)|_{\ve=0}=
p\io |\vv_t|^{p-2}\vv_t\cdot\left(\nabla\xi\cdot\vv_t+\frac{\partial \xi}{\partial t}\right)d\mu_t.$$
Moreover, for $\ve$ sufficiently small (depending on $T$, but not on $t$)
the following inequality holds: 
$$\frac{d}{d\ve}\tilde{V_{\ve}}(t)\leq C (V(\mu,t)+1).$$

\end{lem}
\begin{proof}
If we compute the densities of $\mu^{\ve}_t$ and the expression of the new velocity field and we 
change variable in the integral by $y=T^{\ve}_t(x)$, as we did in the previous lemma, we get
\begin{equation}\label{tildeV}
\tilde{V_{\ve}}(t)=\io \left|\nabla T^{\ve}_t\cdot\vv_t+\frac{\partial T^{\ve}}{\partial t}\right|^p d\mu_t.
\end{equation}
When we differentiate the integrand we get exactly the integrand in \eqref{deriVeq}, and we need only to show
that this expression is uniformly dominated, at least for small $\ve$ and almost every $t$ to get the result.
By boundedness of the derivatives of $T^{\ve}$ it is not difficult to see that the norm of 
the first vector in the scalar product in the integrand may be estimated by
$$\left|\nabla T^{\ve}_t\cdot \vv_t+\frac{\partial T^{\ve}}{\partial t}\right|^{p-1}\leq (C|v_t|+C)^{p-1},$$
while for the second it holds
$$\left|\nabla\xi\cdot \vv_t+\frac{\partial \xi}{\partial t}\right|\leq C|\vv_t|+C$$
for a suitable constant $C$. Hence, since $\vv_t\in [L^p(\mu_t)]^d$ for almost every $t$ 
the integrability is proved and the differentiation under the integral sign can be performed.
\end{proof}

To conclude, we must put together the two previous results in order to 
compute the derivative of the integral in $t$.

\begin{teo}
If $\mu$ is a curve with $\F(\mu)<+\infty$ and $V(\mu,t)\geq V_0>0$ for 
almost every $t$, then it holds
\begin{multline*}
\frac{d}{d\ve}\tilde{\F}(\mu^{\ve})|_{\ve=0}=
\alpha(1-q)\int_0^1F^{\alpha-1}V^{\beta}\io(\nabla\cdot\xi)u_t^q\,d\lcal^d\,dt\\+
p\beta\int_0^1F^{\alpha}V^{\beta-1}\io |\vv_t|^{p-2}(\nabla\xi\cdot \vv_t
+\frac{\partial\xi}{\partial t})\cdot \vv_t \,d\mu_t\,dt,
\end{multline*}
where $F(t)=L_q(\mu_t)$ and $V(t)$ has the usual meaning.
\end{teo}
\begin{proof}
By the definition of $\tilde{\F}(\mu^{\ve})$ we see that the pointwise derivative of the integrand is given by
$\alpha F_{\ve}(t)^{\alpha-1}\frac{dF}{d\ve}\tilde{V_{\ve}}(t)^{\beta}
+\beta F_{\ve}(t)^{\alpha}\tilde{V_{\ve}}(t)^{\beta-1}\frac{d\tilde{V}}{d\ve}$. 
By the regularity of $T^{\ve}$ the term $F_{\ve}(t)$ may be estimated both from above and below by $F(t)$, 
up to multiplicative constants. As far as $\tilde{V}^{\ve}(t)$ is concerned, the argument is a little bit
more tricky. Indeed we must write $\tilde{V}^{\ve}(t)$ according to \eqref{tildeV}, then estimate
$$A^-|\vv_t|-B\leq\left|\nabla T^{\ve}_t\cdot \vv_t+\frac{\partial T^{\ve}}{\partial t}\right|\leq A^+|\vv_t|+B,$$
for $\ve$ small enough,
where the constants $A^\pm$ are as close to $1$ as we want and the constant $B$ is as small as we want 
(this comes from $\nabla T^{\ve}_t=id+O(\ve)$ and $\partial T^{\ve}/\partial t=O(\ve)$), and get
$$A^-\tilde{V}^0-B\leq\tilde{V}^\ve\leq A^+\tilde{V}^{0}+B.$$
The assumption $V\geq V_0>0$ allows us to infer from these inequalities that also 
$\tilde{V}^{\ve}$ may be estimated both from above and below by $V$ up to multiplicative constants. 
Finally, by the estimates in Lemmas \ref{deriF} and \ref{deriV}, we bound the whole pointwise derivative 
by $C F^{\alpha} V^{\beta}$ since we had
$$
\frac{dF}{d\ve}\leq CF;\qquad \frac{d\tilde{V}}{d\ve}\leq C (V +1)
 \leq C(1+\frac{1}{V_0}) V,$$
where the last inequality too comes from $V\geq V_0$. 
Since $L_q^{\alpha} V^{\beta}$ is integrable on $[0,1]$, we may differentiate under the integral
sign and get
$$\frac{d}{d\ve}\F(\mu^{\ve})|_{\ve=0}=
\int_0^1 \left(\alpha F(t)^{\alpha-1}\frac{dF}{d\ve}|_{\ve=0}
\tilde{V}(t)^{\beta}+\beta F(t)^{\alpha}\tilde{V}(t)^{\beta-1}\frac{d\tilde{V}}{d\ve}|_{\ve=0}\right)\,dt.$$
The result follows when we replace the derivatives in $\ve$ by the explicit expressions we computed in 
Lemmas \ref{deriF} and \ref{deriV}.
\end{proof}
\begin{rem}
If $\beta=1/p$ and $\mu$ is a minimizer, it is always possible to get the lower bound $V\geq V_0$ 
by reparametrizing in time, for instance by choosing the constant speed parametrization.
\end{rem}

\begin{cor}
If $\mu$ minimizes $\F$ with given boundary conditions $\mu_0$ and $\mu_1$, 
then its density $u$ and its tangent field $\vv$ satisfy
\begin{multline*}
\alpha(1-q)\int_0^1F(t)^{\alpha-1}V(t)^{\beta}\io(\nabla\cdot\xi)u_t^q\,d\lcal^d\,dt\\
+p\beta\int_0^1F(t)^{\alpha}V(t)^{\beta-1}\io 
u_t|\vv_t|^{p-2}(\nabla\xi\cdot \vv_t+\frac{\partial\xi}{\partial t})\cdot \vv_t \,d\lcal^d\,dt=0,
\end{multline*}
for any vector field $\xi\in C^{\infty}_c(]0,1[\times\Omega;\R^d)$.
\end{cor}
\begin{proof}
It is sufficient to notice that when we create the modified curve
$\mu^{\ve}$ starting form the vector field $\xi$ we do not change the initial and final points of the curve, 
so that the minimality implies that the derivative of $\tilde{\F}(\mu^{\ve})$ at $\ve=0$ vanishes.
\end{proof}

\subsection{A system of PDEs}

The following theorem follows directly from the previous section.

\begin{teo}
Let $\mu_0,\,\mu_1\in W_p(\Omega)$ and let $\mu$ be a curve minimizing $\F$ on $\Gamma(\mu_0,\mu_1)$, with
a finite minimum value.
Then, denoting by $u(t,\cdot)$ the density of $\mu_t$ and by $\vv(t,\cdot)$ the tangent field to the curve $\mu$,
$(u,\vv)$ provide a weak (distributional) solution of the system
\begin{equation}\label{system}
\begin{cases}
        H(t)\nabla u^q+K(t)\nabla\cdot\left(u|\vv|^{p-2}\vv\otimes\vv\right)+
        \frac{d}{dt}\left(K(t)u|\vv|^{p-2}\vv\right)=0& \mbox{ in }\Omega,\\
        \frac{d}{dt}u+\nabla\cdot(\vv u)=0& \mbox{ in }\Omega\\
        u\vv\cdot n=0& \mbox{ on }\partial\Omega\\
        \lim\limits_{t\downarrow 0}u(t,\cdot)\lcal^d=\mu_0;\quad \lim\limits_{t\uparrow 1}u(t,\cdot)\lcal^d=\mu_1,
        \end{cases}
        \end{equation}
where $H(t)=\alpha (1-q)F(t)^{\alpha-1}V(t)^{\beta}$ and $K(t)=p\beta F(t)^{\alpha}V(t)^{\beta-1}$.\\
Given $(\mu_0,\mu_1)$, existence of minimizers is ensured whenever $q<1+1/d$ or,
for general $q$, under the assumption that $\mu_0=u_0\lcal^d$ and $\mu_1=u_1\lcal^d$ with
$u_0,\,u_1\in L^q(\Omega)$ (see \cite{BraButSan}), hence under these conditions existence 
of solutions to this system is ensured.
\end{teo}

It is interesting to rewrite the equations, make some formal simplification and look at some particular cases.

First we expand all the terms in the first equation of System \eqref{system}, obtaining
\begin{multline}\label{svolta}
H(t)\nabla u^q+K(t)\left(u|\vv|^{p-2}v\cdot\nabla\vv+\vv|\vv|^{p-2}\nabla\cdot(u\vv)
+u\left(\vv\cdot\nabla |\vv|^{p-2}\right) \vv \right)\\
+K(t)\left(\vv|\vv|^{p-2}\frac{d}{dt}u
+u\frac{d}{dt}\left(\vv|\vv|^{p-2}\right)\right)+\frac{d}{dt}K(t)u|\vv|^{p-2}\vv=0. 
\end{multline}
Notice that this is always a vector equation, i.e. a system itself, consisting of 
$d$ equations with $d+1$ unknown functions (the components of $\vv$ and the density $u$). 
This system is then completed by the continuity equation. As usual, by $\vv\cdot\nabla \vv$ 
we mean the vector whose $i-$th component is $\sum_j (\vv_j\partial\vv_i/\partial x_j)$.

A formal simplification in \eqref{svolta} may be done: in fact there is a term 
$(K(t)\vv|\vv|^{p-2})(du/dt+\nabla\cdot(u\vv))$ that might be removed by using the continuity equation. 
This is actually possible only under extra regularity assumptions on $K$ and $\vv$ 
(it consists of testing the continuity equation against the product $K(t)\vv|\vv|^{p-2}$ 
which is not in general $C^1$ or regular enough). Anyway, after this formal simplification, 
\eqref{svolta} becomes
\begin{multline}\label{semplificata}
H(t)\nabla u^q+K(t)\left(u|\vv|^{p-2}\vv\cdot\nabla\vv+
u\left(\vv\cdot\nabla |\vv|^{p-2}\right) \vv \right)\\
+K(t)u\frac{d}{dt}\left(\vv|\vv|^{p-2}\right)+\frac{d}{dt}K(t)u|\vv|^{p-2}v=0. 
\end{multline}
Notice that in the case $\beta=1/p$ we can reparametrize in time the solution and there are several 
possible parametrization choices that present some advantages. For instance, we could choose a 
parametrization so that $K(t)$ is constant, to get rid of the final derivative in time. This choice implies
$$V(t)=\left(\frac{F^{\alpha}}{K}\right)^{p/(p-1)},$$
and this, in the case of a bounded $|\Omega|<+\infty$, is sufficient to have the lower bound $V\geq V_0$, 
since in this case $F$ is bounded from below by a positive constant.

Another important fact to be noticed is that in \eqref{semplificata} there is a common $u$ factor. 
It is still formal, but in this way we should get, on $\{u>0\}$,
\begin{eqnarray*}\label{semplificata2}
&&H(t)u^{q-2}\nabla u+K(t)\left(|\vv|^{p-2}\vv\cdot\nabla\vv
+\left(\vv\cdot\nabla |\vv|^{p-2}\right) \vv \right)\\
&&+K(t)\frac{d}{dt}\left(\vv|\vv|^{p-2}\right)+\frac{d}{dt}K(t)|\vv|^{p-2}\vv=0. 
\end{eqnarray*}
\begin{rem}
One might wonder whether the solutions $u$ are automatically positive a.e. in $\Omega$ for $t\in]0,1[$. 
This could be suggested by the fact that in the minimization problem spreadness of the density is favoured. 
In next session we will see with explict examples that this is not necessarily the case.
\end{rem}

We finish this overview of simplifications of the system by looking at the simplest case, i.e. 
$p=q=2$, $\alpha=1$, $\beta=1/2$, in the parametrization regime where $K$ is constant. In this case we get 

\begin{equation}\label{system2}
\begin{cases}
        -2V(t)^{1/2}\nabla u+K\left(\vv\cdot\nabla \vv+\frac{d}{dt}\vv\right)=0& \mbox{ in $\{u>0\}$},\\
        \frac{d}{dt}u+\nabla\cdot(\vv u)=0& \mbox{ in $\Omega$}\\
        u\vv\cdot n=0& \mbox{ on }\partial\Omega\\
        \lim\limits_{t\downarrow 0}u(t,\cdot)\lcal^d=\mu_0;\quad \lim\limits_{t\uparrow 1}u(t,\cdot)\lcal^d=\mu_1.
        \end{cases}
        \end{equation}

Under no constraint on the parametrization we have, instead,
\begin{equation}\label{system3}
\begin{cases}
        -2V(t)^{1/2}\nabla u+K(t)\left(\vv\cdot\nabla\vv+\frac{d}{dt}\vv\right)+\vv\frac{dK}{dt}=0& \mbox{ in $\{u>0\}$},\\
        \frac{d}{dt}u+\nabla\cdot(\vv u)=0& \mbox{ in }\Omega\\
        u\vv\cdot n=0& \mbox{ on }\partial\Omega\\
        \lim\limits_{t\downarrow 0}u(t,\cdot)\lcal^d=\mu_0;\quad \lim\limits_{t\uparrow 1}u(t,\cdot)\lcal^d=\mu_1.
        \end{cases}
        \end{equation}

\section{Self-similar solutions}\label{sec4}

\subsection{Homothetic solutions with fixed center}
In this section we look for particular solutions of the System \eqref{system} 
which are self-similar in the sense that, for any $t$, the measure $\mu_t$ is the image under an homothety 
of a fixed measure.
For simplicity we will consider only the case of System \eqref{system3}, i.e. with $p=q=2$, and we assume that
$0\in\Omega$. The regularity 
of the candidate solutions we will propose will be enough to ensure that we can use this simplified system, 
instead of System \eqref{system}.
To start this analysis it is necessary to establish the following Lemma.
\begin{lem}\label{tangent for ss}
If $\mu$ is a curve in $W_2(\Omega)$ of the form $\mu_t=(T_{R(t)})_{\sharp}\bar{\mu}$ for a 
certain regular function $R:[0,1]\tto ]0,1]$ (where $T_R(x)=Rx$ is the multiplication by a factor $R$, 
hence an homothety), then its tangent field is given by $\vv_t(x)=xR'(t)/R(t)$.
\end{lem}
\begin{proof}
It is not difficult to prove that the field we defined solves the continuity equation and hence is a velocity field.
Indeed, if $\phi\in C^1_c(\Omega)$, it holds
\begin{multline*}
\frac{d}{dt}\io\phi d\mu_t=\frac{d}{dt}\io\phi(R(t)x)\,d \mu(x)
=\io\nabla\phi(R(t)x)\cdot R'(t)x\,d\mu(x)\\=
\io\nabla\phi(R(t)x)\cdot \frac{R'(t)}{R(t)}R(t)x\,d\mu(x)
=\io\nabla\phi\cdot \vv_t \,d\mu_t.
\end{multline*}
It remains to prove that $\vv$ is actually the tangent velocity field, i.e. that its $L^2$ norm is minimal for a.e. $t$. 
This is achieved if we are able to prove that $\Vert\vv_t\Vert_{L^2(\mu_t)}=|\mu|'(t)$ for a.e. $t\in[0,1]$.
To do this, let us fix two times $t<t+h$ and see that the map $T(x)=xR(t+h)/R(t)$ is a transport between 
$\mu_{t}$ and $\mu_{t+h}$. Since it is the gradient of the convex function $x\mapsto x^2R(t+h)/2R(t)$, 
it is actually the optimal transport according to the quadratic cost. Hence
$$
\frac{W_2^2(\mu_t,\mu_{t+h})}{h^2}=\frac{1}{h^2}
\io\left(\frac{R(t+h)}{R(t)}-1\right)^2x^2\,d\mu_t(x)\tto
\io\left(\frac{R'(t)}{R(t)}\right)^2x^2\,d\mu_t(x).
$$
Since this last quantity is exactly the norm of $\vv_t$ in $L^2(\mu_t)$, 
this proves that $\vv$ is the tangent field to the curve $\mu$.
\end{proof}

\begin{rem}
In the case $p\neq 2$ the same result is true, but one has to use the characterization
of tangent velocity fields in terms of closure of gradients of smooth maps, see
Proposition~8.4.5 of \cite{AmGiSa}.
\end{rem}

A first result we prove is the following:
\begin{teo}
If $(u,\vv)$ is a self-similar solution of the system \eqref{system} 
with $u$ Lipschitz continuous, then necessarily $u$ is of the form 
$$u(t,x)=(A_t-B_t|x|^2)\vee 0 \quad\mbox{for suitable coefficients }A_t,\,B_t>0.$$
\end{teo}
\begin{proof}
We look at the equation \eqref{semplificata2} with $p=q=2$, which is valid on $\{u>0\}$, 
and we freeze time, i.e. we look at the resulting space equation for fixed $t$. 
We use the fact that $\vv$ is of the form $\vv_t(x)=c_tx$, which implies that any term $\vv$, $\vv\cdot\nabla\vv$
and $d\vv/dt$ are of the same form. This easily implies that also $\nabla u$ is of the same form. 
Hence, at time $t$, on $\{u>0\}$, it holds $u(x)=A_t-B_tx^2$, where a priori $B_t$ could also be negative. 
Anyway we can prove that $B_t$ cannot be negative. In this case in fact, if $\Omega$ were a convex unbounded 
domain, then $u$ could not be the density of a 
probability measure. On the other hand one can easily see that on bounded convex domains $\Omega$ 
self-similar solutions must vanish on $\partial\Omega$, otherwise we should get a jump of the density 
at the boundary of $\{u>0\}$ when rescaling, but $u$ was supposed to be Lipschitz (except in the case that the solution is constant in time). 
This implies that also in the case of a bounded $\Omega$ the coefficient $B_t$ must be positive. 
For the same continuity reason we get that the region $\{u>0\}$ must agree with the region 
$\Omega\cap\{A_t-B_tx^2>0\}$  in order to have continuity of $u$, and this proves the formula.
\end{proof}

\begin{rem}
A similar result could be obtained for generic Wasserstein spaces with exponent $p>1$, 
getting that any self-similar solution should be of the form $u(t,x)=(A_t-B_t|x|^p)\vee 0$.
\end{rem}

\begin{teo}\label{lorosono}
If $\bar{\mu}$ is a probability measure on $\Omega$ with density $$u(x)=A[(R^2-|x|^2)\vee 0],$$ 
then for any regular and monotone function $R:[0,1]\tto[0,1]$ the curve 
$\mu_t=(T_{R(t)})_{\sharp}\bar{\mu}$ is a solution to the System \eqref{system} 
together with its tangent field $\vv$.
\end{teo}
\begin{proof}
It is sufficient to check the first vector equation in the system \eqref{system3}. 
First we compute the correct constant $A$: we must have
$$1=A\int_0^R(R^2-r^2)d\omega_d r^{d-1}dr=AR^{d+2}\omega_d\frac{2}{d+2},$$
and hence $A=R^{-d-2}(d+2)/(2\omega_d).$ This allows us to compute the term $F(t)$:
$$F=A^2\int_0^R(R^2-r^2)^2d\omega_dr^{d-1}dr=R^{-d}\frac{2(d+2)}{(d+4)\omega_d}.$$
Then we compute $V$ by recalling that $\vv_t(x)=xR'(t)/R(t)$. It holds
$$V=\left(\frac{R'}{R}\right)^2A\int_0^R r^2(R^2-r^2)d\omega_d r^{d-1}dr=\frac{d}{d+4}(R')^2.$$
We must also compute $d\vv/dt$ and $\vv\cdot\nabla \vv$:
$$\frac{\partial \vv}{\partial t}=x\frac{R''R-(R')^2}{R^2};
\quad \nabla\vv=\left(\frac{R'}{R}\right)I;\quad \vv\cdot\nabla \vv=\left(\frac{R'}{R}\right)^2x.$$
We compute now 
\begin{gather*}
K(t)=F(t)V(t)^{-1/2}=R^{-d}|R'|^{-1}\frac{2(d+2)}{\sqrt{d(d+4)}\omega_d},\\
K'(t)=\sign(R')(-dR^{-d-1}-R^{-d}(R')^{-2}R'')\frac{2(d+2)}{\sqrt{d(d+4)}\omega_d}.
\end{gather*}
If we call $c=\sign(R')\frac{2(d+2)}{\sqrt{d(d+4)}\omega_d}$ it holds 
$K=cR^{-d}(R')^{-1}$ and $K'=c(-dR^{-d-1}-R^{-d}(R')^{-2}R'')$, 
but also $-2V^{1/2}\nabla u(x)=cdR'R^{-d-2}x$.
Inserting everything in the equation we must check that
$$
dR'xR^{-d-2}+R^{-d}(R')^{-1}x\frac{R''}{R}-(dR^{-d-1}+R^{-d}(R')^{-2}R'')x\frac{R'}{R}=0.$$
The proof is achieved as this last equation is (miracolously enough)
always satisfied.\end{proof}

\begin{rem}
By a similar proof we can show that, for $p\neq 2$, if $\bar{\mu}$ has a density of the form 
$u(x)=A[(R^p-|x|^p)\vee 0]$, then $\mu$ gives raise to a self-similar solution.
\end{rem}

\begin{rem}
This kind of self-similar solutions can join two different probability measures which are homothetic, and 
in particular arrive up to the Dirac mass $\delta_0$. Anyway it is not in general possible to link a measure 
to $\delta_0$ by a curve with finite energy: in \cite{BraButSan}, conditions to ensure this possibility 
are provided, but in general they are not satisfied in the case $q=2$. 
\end{rem}

\subsection{Moving self-similar solutions}
We have characterized all the self-similar solutions which link two homothetic probability measures. 
It is however interesting to look also at the moving self-similar solutions, i.e. at solutions 
obtained by homotheties and translations together.

In this case we consider a reference measure $\bar{\mu}$ and we look for solutions of the form 
$(T^t)_{\sharp}\bar{\mu}$, where $T^t(x)=R(t)x+\bar{x}(t)$. 
It is not difficult to replace Lemma \ref{tangent for ss} with the following:
\begin{lem}
If $\mu$ is a curve of the form $\mu_t=(T^t)_{\sharp}\bar{\mu}$, then its tangent field is given by 
$$\vv_t(x)=\frac{R'(t)}{R(t)}(x-\bar{x}(t))+\bar{x}'(t).$$
\end{lem}
\begin{proof}
The result may be proved very similarly to Lemma \ref{tangent for ss}: 
it is sufficient to check the continuity equation
\begin{multline*}
\frac{d}{dt}\io\phi(R(t)x+\bar{x}(t)) \,d\mu(x)
=\io\nabla\phi(R(t)x+\bar{x}(t))\cdot (R'(t)x+\bar{x}'(t))\,d\mu(x)\\=
\io\nabla\phi(R(t)x+\bar{x}(t))\cdot \frac{R'(t)}{R(t)}(R(t)x+\bar{x}'(t))\,d\mu(x)
=\io\nabla\phi\cdot \vv_t \,d\mu_t,
\end{multline*}
and then to check the optimality of the norm by the fact that the map 
$$x\mapsto \frac{R(t+h)}{R(t)}(x-\bar{x}(t))+\bar{x}(t+h)$$
transports $\mu_t$ on $\mu_{t+h}$ and is optimal, and that
$$\frac{1}{h^2}\io\left(\frac{R(t+h)}{R(t)}(x-\bar{x}(t))+\bar{x}(t+h)-x\right)^2\,d\mu_t(x)$$
converges to
$$
\io\left(\frac{R'(t)}{R(t)}
(x-\bar{x}(t))+\bar{x}'(t)\right)^2\,d\mu_t(x)=\Vert\vv_t\Vert_{L^2(\mu_t)}.
\qedhere
$$
\end{proof}
For computational simplicity we consider moving self-similar solutions only under a special
reparametrization.
\begin{teo}
If $\bar{\mu}$ is a probability measure on $\Omega$ with density $$u(x)=A[(R^2-|x|^2)\vee 0]$$ 
and $\bar{x}(0),\,\bar{x}(1)\in\Omega$ are assigned, a curve $\mu_t=(T^t)_{\sharp}\bar{\mu}$, 
parametrised so that $K=FV^{-1/2}$ is constant, is a moving self-similar solution 
(solving System \eqref{system2} together with its own tangent field) if and only if 
the vector $x$ moves on the straight line segment from $\bar{x}(0)$ to $\bar{x}(1)$ 
with constant speed and $R$ is a strictly concave function of $t$. This means
$$\bar{x}''=0;\quad R^{2d}(d(R')^2+(d+4)(\bar{x}')^2)\,\mbox{ is constant   and $R$ strictly concave.}$$
\end{teo}
\begin{proof}
We only need to check under which conditions the first equation is satisfied. 
We re-write in this case the quantity considered in Theorem \ref{lorosono}: first we compute
\begin{gather*}
u(x)=A[(R^2-|x-\bar{x}|^2)\vee 0];\quad A=\frac{(d+2)}{2R^{d+2}\omega_d};\quad 
\nabla u(x)=-\frac{(d+2)}{R^{d+2}\omega_d}(x-\bar{x})\\
F=R^{-d}\frac{2(d+2)}{(d+4)\omega_d};\quad V=\frac{d}{d+4}(R')^2+(\bar{x}')^2.
\end{gather*}
We have used the fact that $u_t$ is symmetric around $\bar{x}(t)$ and hence there is no 
mixed term $(x-\bar{x}(t))\cdot\bar{x}'(t)$ in computing $V(t)$. 
Then we go on with $d\vv/dt$ and $\vv\cdot\nabla\vv$:
\begin{gather*}
\frac{\partial \vv}{\partial t}
=(x-\bar{x})\frac{R''R-(R')^2}{R^2}-\bar{x}'\frac{R'}{R}+\bar{x}'';\quad 
\nabla \vv=\left(\frac{R'}{R}\right)I;\\
\vv\cdot\nabla \vv=\left(\frac{R'}{R}\right)^2(x-\bar{x})+\frac{R'}{R}\bar{x}'\quad
\frac{\partial \vv}{\partial t}+\vv\cdot\nabla\vv=(x-\bar{x})\frac{R''}{R}+\bar{x}''.
\end{gather*}
Then we look at the the condition to have $K'(t)=0$, which is equivalent to $F^{-2}V$ being constant, 
and thus $R^{2d}(d(R')^2+(d+4)(\bar{x}')^2)$ must be constant.
Assuming $K$ to be constant we try to satisfy the equation, and we write it in the following form that 
we can reach after multiplying by $V^{1/2}$:
$$-2V\nabla u+F\left(\frac{\partial \vv}{\partial t}+\frac{1}{2}\vv\cdot\nabla\vv\right)=0.$$
This equation becomes
$$2\left(\frac{d}{d+4}(R')^2+(\bar{x}')^2\right)\frac{(d+2)}{\omega_dR^{d+2}}(x-\bar{x}(t))
+R^{-d}\frac{2(d+2)}{(d+4)\omega_d}((x-\bar{x})\frac{R''}{R}+\bar{x}'')=0.$$
To satisfy this equation it is necessary and sufficient that the two parts, 
the one involving $x-\bar{x}$ and the other with $\bar{x}''$ both vanish. After simplifying we get
$$R^{-2}(d(R')^2+(d+4)(\bar{x}')^2)+\frac{R''}{R}=0;\qquad \bar{x}''=0.$$
Hence we must have $\bar{x}(t)=(1-t)\bar{x}(0)+t\bar{x}(1)$ and $\bar{x}'(t)=e=\bar{x}(1)-\bar{x}(0)$.
Now we recall that $R^{2d}(d(R')^2+(d+4)(\bar{x}')^2)$ was assumed to be constant and so 
$d(R')^2+(d+4)(\bar{x}')^2=CR^{-2d}$. Hence we get $R''=-CR^{-2d-1}$. 
Thus, $u$ is a moving self-similar solutions if and only if the following conditions simultaneously hold:
$$\begin{cases}d(R')^2+(d+4)e^2=CR^{-2d}&\mbox{for a certain } C,\\
               R''=-CR^{-2d-1}&\mbox{for the same } C,\\
               \bar{x}(t)=\bar{x}(0)+te.&\end{cases}$$
By differentiating the first equation we get $2dR'R''=-2dCR^{-2d-1}R'$ and hence the second is automatically 
satisfied, provided we can ensure that $R'\neq 0$ a.e. This means that $R$ being strict concave is sufficient 
(it is not possible to have more than a time where $R'$ vanishes), 
but it is also necessary from the second equation. The result is then proved.
\end{proof}


\begin{thebibliography}{999}

\bibitem{AmGiSa} L. Ambrosio, N. Gigli and G. Savar\'e,
{\it Gradient flows in metric spaces and in the spaces of probability measures.} 
Lectures in Mathematics, ETH Zurich, Birkh\"auser, 2005.

\bibitem{BraButSan} A. Brancolini, G. Buttazzo and F. Santambrogio, Path Functionals over Wasserstein 
spaces, 2005. {\it J. Eur. Math. Soc.}, to appear, available at {\tt cvgmt.sns.it}.

\bibitem{Least} Y. Brenier, The Least Action Principle and the Related Concept of Generalized Flows 
for Incompressible Perfect Fluids. {\it J. Amer. Math. Soc.}, {\bf 2} (1989), 225--255.

\bibitem{Br1} Y. Brenier, The dual least action principle for an ideal incompressible fluid.
{\it Arch. Rational Mech. Anal.}, {\bf 122} (1993), 323--351.

\bibitem{Br2} Y. Brenier, A homogenized model for vortex sheets.
{\it Arch. Rational Mech. Anal.}, {\bf 138} (1997), 319--353.

\bibitem{Br3} Y. Brenier, Minimal geodesics on groups of volume-preserving maps and generalized solutions
of the Euler equation. {\it Comm. Pure and Appl. Math.}, {\bf 52} (1999), 411--452.

\bibitem{guionnet} A. Guionnet, First order asymptotics of matrix integrals; a rigorous approach
towards the understanding of matrix models. {\it Comm. Math. Phys.}, {\bf 244} (2004), 527--569.

\bibitem{loeper} G. Loeper, The reconstruction problem for the Euler-Poisson system in cosmology.
{\it Arch. Rational Mech. Anal.}, to appear, available at {\tt http://math.univ-lyon1.fr/$\sim$loeper}.

\bibitem{MadMorSol} F. Maddalena, S. Solimini and J.-M. Morel,  A variational model of irrigation patterns. {\it Interfaces and Free Boundaries} {\bf 5} (2003), 391--415. 

\bibitem{Villani} C. Villani, Topics in optimal transportation.
{\it Graduate Studies in Mathematics} {\bf 58},
American Mathematical Society, 2003.

\bibitem{xia1} Q. Xia, Optimal Paths related to Transport Problems. {\it Comm. Cont. Math.} {\bf 5} (2003), no. 2, 251--279.

\end{thebibliography}
\end{document}